\documentclass[a4paper,11pt]{amsart}
\usepackage[lmargin=2.5cm,tmargin=2.5cm,bmargin=2.5cm,rmargin=2.5cm]{geometry}
\usepackage[utf8]{inputenc}
\usepackage{indentfirst, amsfonts, amsmath, amsthm, amssymb, amscd,mathtools,amscd,bezier}
\usepackage[dvips]{graphicx}
\usepackage{xcolor}
\usepackage{wrapfig}
\usepackage{pstool}
\usepackage{overpic}
\usepackage{MnSymbol}
\usepackage{graphics}
\usepackage{epsfig}
\usepackage{psfrag}
\usepackage{hyperref}
\usepackage[normalem]{ulem}
\usepackage{thmtools}
\usepackage{thm-restate}
\usepackage{setspace}

\newtheorem{maintheorem}{Theorem}

\newtheorem{theorem}{Theorem}[section]

\newtheorem{corollary}{Corollary}
\newtheorem{proposition}{Proposition}

\newtheorem{definition}[theorem]{Definition}

\theoremstyle{remark}
\newtheorem{remark}[theorem]{Remark}

\newcommand{\R}{\mathbb{R}}

\newcommand{\N}{\mathbb{N}}

\renewcommand{\S}{\mathbb{S}}
\newcommand{\s}{\Sigma}
\newcommand{\g}{\Gamma}

\title[There exist transitive PSVF on $\mathbb{S}^2$ but not robustly transitive.]{There exist transitive piecewise smooth vector fields on $\mathbb{S}^2$ but not robustly transitive.}

 \author{Rodrigo D. Euz\'ebio}
 \address{Instituto de Matemática e Estatística, IME-UFG, Goi\^{a}nia-GO, Brazil.}
   \email{euzebio@ufg.br}
   \author{Joaby S. Juc\'{a}}
   \address{Instituto de Matemática e Estatística, IME-UFG, Goi\^{a}nia-GO, Brazil.}
   \email{joabyjuca@discente.ufg.com}
 \author{R\'egis Var\~{a}o} 
 \address{Departamento de Matem\'atica, Estat\'istica e Computa\c c\~ao Cient\'ifica,
   IMECC-UNICAMP, Campinas-SP, Brazil.}
 \email{varao@unicamp.br}

\begin{document}

\doublespacing

 \maketitle

\begin{abstract} 
It is well known that smooth (or continuous) vector fields cannot be topologically transitive on the sphere $\S^2$. Piecewise-smooth vector fields, on the other hand, may present non-trivial recurrence even on $\S^2$. Accordingly, in this paper the existence of topologically transitive piecewise-smooth vector fields on $\S^2$ is proved, see Theorem \ref{teorema-principal}. We also prove that transitivity occurs alongside the presence of some particular portions of the phase portrait known as {\it sliding region} and {\it escaping region}. More precisely, Theorem \ref{main:transitivity} states that, under the presence of transitivity, trajectories must interchange between sliding and escaping regions through tangency points. In addition, we prove that every transitive piecewise-smooth vector field is neither robustly transitive nor structural stable on $\S^2$, see Theorem \ref{main:no-transitive}. We finish the paper proving Theorem \ref{main:general} addressing non-robustness on general compact two-dimensional manifolds.

\end{abstract}

\vspace{.5cm}


\smallskip


\section{Introduction}\label{setting-the-problem}

Piecewise-smooth vector fields (PSVF) describe a particular type of dynamical systems for which discontinuities may occur on the phase portrait. Through the literature several authors have dealt with PSVF by assuming distinct conventions in order to establish how trajectories interact with those discontinuities. A particular approach largely adopted in the study of PSVF was established by A. F. Filippov in \cite{Fi}. The convention adopted by Filippov captures the so-called {\it sliding motion} in such way that trajectories may reach the discontinuity set in finite time and then slide on that region. An equivalent approach is considered by Utkin (see \cite{Utkin1992}). Other conventions are possibly more restrictive in the treatment of the trajectories, we mention Broucke et. al. (see \cite{BrouckePughSimic01}) and Barbashin (see \cite{Barbashin1970}). In this paper we adopt Filippov's convention. As a motivation to the study of PSVF on the background of Filippov, we mention that several practical problems can be modeled by such vector fields, for instance stick-slip processes, the anti-lock braking system (ABS), the relay systems and generally speaking a substantial part of control systems (see these and other applications of PSVF in \cite{diBernardo-electrical-systems}, \cite{diBern-relay}, \cite{Brogliato}, \cite{Rossa}, \cite{Dixon}, \cite{Genema},  \cite{Jac-To}, \cite{Kousaka} and \cite{Leine}).

 
An important feature of smooth and piecewise-smooth vector fields concerns topological transitivity. The classical literature on dynamical systems has a well established theory on transitivity for smooth vector fields (an, of course, for diffeomorphims), but such a theory is still premature in the particular context of PSVF. Indeed, some preliminary results can be found in  \cite{BCE,BCE-ETDS,Carvalho-LFernando, EV} and references therein. A PSVF to be topologically transitive requires that any two open sets can be connected by a Filippov orbit (see the proper definitions on Section \ref{sec:preliminaries}).

It is a well known result that the two-dimensional sphere $\mathbb{S}^2$ does not admit topologically transitive continuous vector fields (see \cite{Lopez-Lopez04} and references therein). In this paper, we investigate that and other related issues for PSVF defined on the sphere. In particular, we are able to provide an explicit example of an one-parameter family of topologically transitive PSVF on $\mathbb{S}^2$.

\begin{maintheorem}\label{teorema-principal}
	There exist an one-parameter family of transitive piecewise-linear vector fields defined on the sphere $\S^2$.
\end{maintheorem}

\begin{figure}
	\includegraphics[scale=0.4]{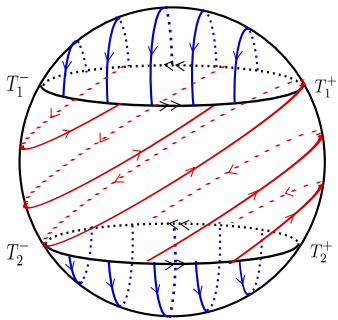}
	\caption{Filippov orbits of a topologically transitive PSVF on $\S^2$}
	\label{main-figure}
\end{figure}

One may notice from Figure \ref{main-figure} that the topologically transitive vector fields we provide in Theorem \ref{teorema-principal} present three zones separated by two circles. A proper question concerns the possibility of constructing a transitive PSVF on the two-dimensional sphere with only two zones. Actually, we are able to partially answer this question. The following result states that, for a certain class of PSVF, the answer is negative when considering two zones separated by a circle.


\begin{proposition}\label{linear-2zones}
There exist no transitive piecewise-linear vector fields on $\S^2$ having two zones and separated by a circle.
\end{proposition}

Transitive PSVF on $\S^2$ have some inherent features concerning the discontinuities occurring in the phase portrait. We point to the existence of sliding and escaping regions (see Section \ref{sec:preliminaries} for precise definitions) connecting each other through tangency points.

 \begin{maintheorem}\label{main:transitivity}
 If $Z$ is a transitive PSVF  on $\S^2$ having a finite number of tangency points on $\Sigma$, then the following statements hold:
	\begin{itemize}
		\item[(a)] the sliding and escaping regions are non-empty sets;
		\item[(b)] every sliding and escaping regions are connected by some trajectory of $Z$. Moreover there are an uncountable number of trajectories of $Z$ connecting sliding and escaping regions.
	\end{itemize}
\end{maintheorem}

Theorems \ref{teorema-principal} and \ref{main:transitivity} together illustrate the richness of the trajectories in a PSVF. Nevertheless, every transitive PSVF on $\S^2$ must present a non-trivial recurrence between two type of sets which arise when discontinuities are allowed in the vector fields. On the other hand, such a recurrence is easily broken by small perturbations, that is, transitivity is not a generic property for PSVF as stated in the next result.


\begin{maintheorem}\label{main:no-transitive}
 There exist no robustly transitive PSVF on $\S^2$ with finite number of tangency points.
\end{maintheorem}


We are mainly concerned with $\mathbb S^2$ but our techniques can be effortlessly applied to obtain results for compact two-dimensional manifolds. 

 \begin{maintheorem}\label{main:general}
 Let $M^2$ be a two-dimensional compact manifold. There exist no robustly transitive PSVF on $M^2$ having non-empty sliding and escaping regions with finite number of tangency points.
 \end{maintheorem}

The next result is a direct consequence of the non-robustness of transitive PSVF.

\begin{corollary}\label{main:no-structural-stable}
Every transitive Filippov PSVF defined on a two-dimensional compact manifold with finite number of tangency points is not structurally unstable.
\end{corollary}

The corollary is proved by noticing that in order to a  PSVF fitting the hypotheses of Theorem \ref{main:no-transitive} be structurally stable, it should also be robustly transitive which is a contradiction to Theorem \ref{main:no-transitive}. For the sake of completeness we remark that a different proof of Corollary \ref{main:no-structural-stable} can be achieved using the results on structural stability presented in \cite{BrouckePughSimic01}, which applies to Filippov convention due to a suitable definition of topological conjugation considered in that paper. For our purposes only the presented proof will be enough.

The next section contains the precise definitions used throughout this paper and the last one contains the proofs of the main results.

\section{Preliminaries}\label{sec:preliminaries}

\begin{definition}\label{definicao_PSVF}
    A piecewise-smooth vector field is a triple $(M,\Sigma,Z)$ where
    \begin{itemize}
	\item[(i)] $M$ is a suitable manifold;
	\item[(ii)] $\Sigma$ is formed by a finite union of simple curves $\s=\s_1\cupdot\cdots\cupdot\s_n$ splitting $M$ into $n+1$ connected components regions $R_i$, where $\s_i=\gamma_i^{-1}(0)$ and $\gamma_i:M\rightarrow \mathbb{R}$ are smooth functions having $0$ as regular value, $i=1,\ldots,n$;
	\item[(iii)] $Z$ is a collection of $n+1$ vector fields of class $C^r$ defined on $M$, say $Z=(X_1,\ldots,X_{n+1})$, being each $X_i$ defined on the closure of $R_i$. 
\end{itemize}
\end{definition}

We shall denote a PSVF by $Z$ instead of the triple  $(M,\Sigma,Z)$ unless there is some confusion on $M$ or $\Sigma$. We call $\Sigma$ the \textbf{switching manifold} and we notice that $Z$ is bi-valuated on $\Sigma$. In particular, every component $X_i$ of $Z$ is a vector field defined on whole $M$ which has been restricted to $R_i$. Because $Z$ is bi-valuated on each connected component of $\Sigma$, it is necessary to establish some {\it rule} describing how trajectories interact to $\Sigma$, switching to one side of $\Sigma$ to another or even remaining on it. As mentioned before, in this paper we adopt the Filippov convention which we describe in this section.

\begin{remark}
We remark that Filippov convention requires connected components of $\Sigma$ to be disjoint pairwise, simple and smooth. Nevertheless other conventions or including extensions of the Filippov one may suppress some of those assumptions in such way that Definition \ref{definicao_PSVF} could be slightly adapted.
\end{remark}




Let $\s_i= \gamma_i^{-1}(0)$ be the common boundary between the regions $R_i$ and $R_j$ and suppose that $\nabla \gamma_i(p)$ points to the interior of the region $R_i$ for all $p\in\s_i$. 
We distinguish three regions on $\s_i$ satisfying $(X_i.\gamma_i(p))\cdot (X_j.\gamma_i(p))\neq 0$, where $X_k.\gamma_i(p)=\langle X_k(p),\nabla \gamma_i(p)\rangle$ is the first Lie derivative of $\gamma_i$ in the direction of vector field $X_k$ at the point $p$. Such regions are characterized in what follows.

\begin{definition}
Let $\Sigma_i$ be a connected component of $\Sigma$ for some $i=1,\ldots,n$. Call $X_i$ and $X_j$ the vector fields separated by $\Sigma_i$ and let $p\in\Sigma_i$ such that $X_i(p)$ and $X_j(p)$ are transversal to $\Sigma_i$ at $p$.  Under this assumptions we distinguish three types of regions on $\s_i$:
    \begin{itemize}
        \item[i)] The {\bf crossing region $\Sigma_i^c$ of $\Sigma_i$} which is formed by the points $p\in\Sigma_i$ such that $(X_i.\gamma_i(p))\cdot (X_j.\gamma_i(p))>0$. 
        \item[ii)] The {\bf escaping region $\Sigma_i^e$ of $\Sigma_i$} which is formed by the points $p\in\Sigma_i$ such that $X_i.\gamma_i(p)>0$ and $X_j.\gamma_i(p)<0$.
        \item[iii)] The {\bf sliding region $\Sigma_i^s$ of $\Sigma_i$} which is formed by the points $p\in\Sigma_i$ such that $X_i.\gamma_i(p)<0$ and $X_j.\gamma_i(p)>0$.
    \end{itemize}
We call $\Sigma^c=\displaystyle\bigcup_{i=1}^{n}\Sigma_i^c$, $\Sigma^s=\displaystyle\bigcup_{i=1}^{n}\Sigma_i^s$ and $\Sigma^e=\displaystyle\bigcup_{i=1}^{n}\Sigma_i^e$ the crossing, sliding and escaping regions of $\Sigma$, respectively.
\end{definition}

Notice that, when $\nabla \gamma_i(p)$ points to the interior of $R_j$ for all $p\in\s_i$ the inequalities in bullets $ii)$ and $iii)$ of the last definition are interchanged.


\begin{definition}
The points $p\in\Sigma_i$ such that $X_i. \gamma_i(p) = 0$ (resp. $X_j. \gamma_i(p) = 0$) are called \textit{tangency points} of $X_i$ (resp. $X_j$). The collection of the points $p\in\Sigma$ such that $p$ is a tangency point for some vector field $X_j$, $j=1,\ldots,n+1$ constitute the {\bf set of tangency points} of $Z$ denoted by $\Sigma^t$.
\end{definition}


 Let $p\in\Sigma$ be a tangency point of $Z=(X_1,\ldots,X_{n+1})$. We say that $Z$ has a contact of order $n\in\N$ with $\Sigma$ at $p$ if, for some $X_i$, $X_i^k.\gamma_i(p)=\langle \nabla X_{i}^{k-1}.\gamma_i(p),X_i(p)\rangle=0$ for $k<n$ and $X_i^n.\gamma_i(p)\neq 0$. We classify tangency points according to the following: We say that $p\in\Sigma$ is an {\it invisible tangency point} if $X_i$ has a contact of even order at $p$ and $X_i^r.\gamma_i(p)<0$. On the other hand, we say that $p\in\Sigma$ is a {\it visible tangency point} if either $X_i$ has a contact of even order at $p$ and $X_i^r.\gamma_i(p)>0$ or $X_i$ has contact of odd order at $p$.
 
 A particular kind of tangency points are the even contact order points which are tangency points for both $X_i$ and $X_j$. We refer to those points by {\it double tangency}. We say that a double tangency $p$ is {\it elliptic} if $p$ is invisible for $X_i$ and $X_j$, {\it hyperbolic} if it is visible for $X_i$ and $X_j$ and {\it parabolic} if $p$ is visible for $X_i$ and invisible for $X_j$ or otherwise.

\begin{figure}
\includegraphics[scale=1.5]{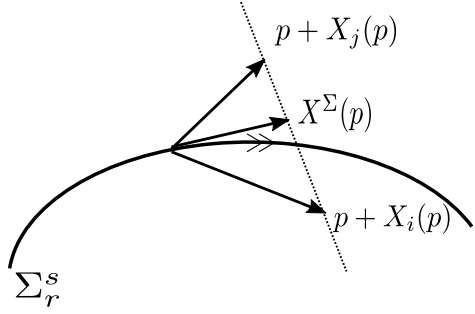}
\caption{Filippov's convention}
\label{filippov-convention}
\end{figure}

In what follows we define the sliding vector field $Z^\s$ on $\s^{s\cup e}$. If $\s_i$ is a common boundary of $R_i$ and $R_j$ and $p\in\s_i^s$ then we define $Z^\s(p)=m-p$, with $m$ being the point on segment joining $p+X_i(p)$ and $p+X_j(p)$ such that $m-p$ is tangent to $\s_i^s$ (see Figure \ref{filippov-convention}). If $p\in\s_i^e$, then $p\in\s_i^s$ for the vector field $-X$ and we define $Z^\s(p)=-(-X)^\s(p)$. In our pictures we represent the dynamics of $Z^\s$ by double arrows and we refer to $Z^\s$ by {\it Filippov vector field}. The points $p\in\Sigma^{s,e}$ such that $Z^\Sigma(p)=0$, that is, the equilibrium points of the Filippov vector field, are called \textit{pseudo equilibrium points} of $Z$.

\begin{remark}
	We notice that although the Filippov vector field is defined at sliding and escaping points we can extend it beyond the boundary of $\s^{s,e}$. For instance, if $p\in\Sigma^t$ and
	$$
	\lim_{q\to p}X^\s(q)=L\neq 0,\ \ q\in\s^{s,e},
	$$
	then we define the extended Filippov vector field at $p$ as $X^\s(p)=L$. That will be the case in the proof of Theorem \ref{teorema-principal}.
\end{remark}

A \textit{global trajectory} $\Gamma_Z(t,p)$ of a piecewise-smooth vector field $Z$ is the trace of a 	continuous curve obtained by a suitable oriented concatenation of trajectories of $X_i$ and/or $X_j$ and/or $Z^\Sigma$. A \textit{maximal trajectory} $\Gamma_Z(t,p)$ is a global trajectory that cannot be extended by any concatenation of trajectories of $X_i,\ X_j$ or $Z^\Sigma$. We also refer to a maximal trajectory $\Gamma_Z(t,p)$ as a {\it Filippov trajectory (or orbit)}. If a Filippov trajectory is singular at $p\in M$ we say that $p$ is an equilibrium point of $Z$ in the sense that $X_i(p)=0$ for some $X_i$ which is defined on $R_i$. We say that $p$ is a {\it real} equilibrium point if $p$ belongs to the closure of $R_i$, otherwise we say that $p$ is a {\it virtual} equilibrium point.


We finish this section by introducing two definitions addressing topolocally transitive PSVF which are inspired in the classical definitions for smooth vector fields.

\begin{definition}
	A PSVF is \textit{topologically transitive} if given two arbitrary open sets $U$ and $V$ of $M$, there exist a Filippov trajectory connecting these sets.
\end{definition}



\begin{definition}
We say that a $Z$ is robustly topologically transitive if $Z$ is topologically transitive and every PSVF sufficiently close to $Z$ is also topologically transitive.

\end{definition}


\section{Proof of the main results}\label{proofs}

In this section we prove the main results of the paper. We consider the unitary sphere $\S^2\subset\R^3$ centered at the origin.

\subsection{Proof of Theorem \ref{teorema-principal}}
Initially we construct a topologically transitive pie\-ce\-wi\-se-linear vector field on $\S^2$ (see Figure \ref{main-figure}). Then we perturb it to obtain the desired one-parametric family of transitive PSVF on $\S^2$.
	
Effectively, let $\s_1$ and $\s_2$ be the curves on $\S^2$ given by the intersection of $\S^2$ with the planes $z=1/2$ and $z=-1/2$, respectively. Consider the linear vector fields
$$
X(p)=(z,0,-x)\quad\mbox{and}\quad Y(p)=\left(-\frac{1}{2}\left(\sqrt{3}\,y+z\right),\frac{\sqrt{3}}{2} x, \frac{x}{2}\right)
$$
with $p=(x,y,z)\in\S^2$. Let $Z=(X,Y,X)$ be a piecewise-linear vector field with three zones on $\S^2$ being
$X$ defined on $R_1=\{(x,y,z)\in\S^2;\ z\geq1/2\}$ and $R_3=\{(x,y,z)\in\S^2;\ z\leq-1/2\}$ and $Y$ defined on $R_2=\{(x,y,z)\in\S^2;\ |z|\leq1/2\}$. Notice that $Y$ is the vector field obtained from $-X$ trough the rotation by the angle $\pi/3$ around of the $x-$axis in the clockwise sense. We claim that $Z$ is topologically transitive on $\S^2$.
	
 In order to prove the claim we first check that $Z$ satisfies the following properties:
	\begin{itemize}
		\item[(i)] Each equilibrium point of $Z$ is virtual;
		\item[(ii)] $Z$ has a pair of double tangency points $T_1^-\neq T_1^+$ with $\s_1$, a pair of double tangency points $T_2^-=-T_1^+$ and $ T_2^+=-T_1^-$ with $\s_2$, and $Z$ has no more tangency points with $\s$ besides $T_i^\pm,\ i=1,2$;
		\item[(iii)] the tangency points $T_i^\pm,\ i=1,2$ are invisible for $X$, $T_1^+$ and $T_2^-$ are visible for $Y$ and $T_1^-$ and $T_2^+$ are invisible for $Y$;
		\item[(iv)] $Y$ has a periodic orbit on $R_2\cup\s$ connecting $T_2^-$ to $T_1^+$;
		\item[(v)] $\s^c=\emptyset$; $\s_1^e=\{(x,y,z)\in\s_1;\ x<0\}$, $\s_1^s=\{(x,y,z)\in\s_1;\ x>0\}$, $\s_2^e=\{(x,y,z)\in\s_1;\ x<0\}$ and $\s_2^s=\{(x,y,z)\in\s_1;\ x>0\}$;
		\item[(vi)] the extended Filippov vector field is well defined on the entire $\s_{1,2}$  and its orientation is in the counter-clockwise sense without pseudo-equilibrium points.
	\end{itemize}
	
	The claims $(i)-(v)$ can be proved straightforward by a direct integration of the linear vector fields $X$ and $Y$. To prove claim $(vi)$, we notice that $X_t(p) = (x \cos(t) + z \sin(t),y,z \cos(t) - x \sin(t))$ and $Y_t(p) = (x_t(p),\ y_t(p),\ z_t(p))$ where
	$$
	\begin{array}{lcll}
	x_t(p) = x \cos(t) - 1/2 (\sqrt{3} \,y + z) \sin(t); \vspace{0.2cm}\\
	y_t(p) = \frac{1}{4} (y + 3 y \cos(t) + \sqrt{3} (z (-1 + \cos(t)) + 2 x \sin(t))); \vspace{0.2cm}\\
	z_t(p) = \frac{1}{4} (3 z + \sqrt{3}\, y (-1 + \cos(t)) + z \cos(t) + 2 x \sin(t))),
	\end{array}
	$$
	with $p=(x,y,z)\in\S^2$. Consequently, the extended Filippov vector field on $\s_i,\ i=1,2$, can be defined by
	\begin{equation}\label{campo-extendido-de-Filippov}
	Z^\s(p)=\dfrac{Yf(p)X(p)-Xf(p)Y(p)}{Yf(p)-Xf(p)}=\dfrac{\sqrt{3}}{3}(-y,x,0),
	\end{equation}
	for $p=(x,y,z)\in\s_i^{e,s},\ i=1,2,$ and therefore claim $(vi)$ is proved.
	
	Now we prove the transitive property of $Z$. Notice that, for every $p=(x,y,z)$ with $|z|\geq\frac{1}{2}$, $p\neq T_i^\pm,\ i=1,2$, there exist $t>0$ such that $X_t(p)\in\s_1^s$. Analogously, for every $p=(x,y,z)$ with $|z|\leq\frac{1}{2}$, $p\neq T_i^\pm,\ i=1,2$, there exist $t>0$ such that $Y_t(p)\in\s_2^s$. On the other hand, if $p\in\s$, then there exist $t>0$ such that either $Z^\s_t(p)=T_1^+$ for $p\in\s_1$ or $Z^\s_t(p)=T_2^-$ for $p\in\s_2$.
	
	Moreover, there exist a periodic orbit for $Z$ on the region $R_2$ connecting $T_1^+$ to $T_2^-$ (see statement (iv)), then given $p,q\in\S^2$ there exists $t_1,\ t_2>0$ and a maximal trajectory $\g_Z(t,p)$ satisfying $\g_Z(t_1,p) = q = \g_Z(-t_2,p)$. Therefore, given two arbitrary open set $U$ and $V$ on $S^2$, we can connect every point on $U$ to every point on $V$ so $Z$ is transitive in $\S^2$.

From $Z$ we construct an one-parametric family $Z_\theta=(X,Y_\theta,X)$ by setting on the central region the family of vector fields $Y_\theta$ with $\pi/6<\theta\leq\pi/3$. This new vector field is obtained by the rotation of $-X$ by an angle $\theta$ in the clockwise sense around the $x$-axis. The linear center of $Y_\theta$ remains virtual and therefore the new PSVF $Z_\theta$ on $\S^2$ is transitive for $\pi/6<\theta\leq\pi/3$.

For $\theta\leq\pi/6$ the center lies in $\s$ or it is a real equilibrium of $Z$ and for $\theta>\pi/3$ there exist no trajectory of $R_i$ visiting $R_j$ for $i,j=1,2,\ i\neq j$ so no longer transitive regime takes place.


\subsection{Proof of Proposition \ref{linear-2zones}} Let $Z=(X_1,X_2)$ be a linear PSVF defined on $\S^2$ separated by a single circle $\Sigma$. Before consider the problem of transitivity, we notice that once $X_i(\bf x)=A_i\bf x$ is linear, using straightforward calculations from linear algebra, $A_i$ is a real skew-symmetric matrix so the eigenvalues of $A_i$ are zero and a pair of conjugated imaginary pure eigenvalues, $i=1,2$. Thus $X_1$ has an one-dimensional invariant straight line trough the origin filled by equilibrium points associated to the 0 eigenvalue of $A_1$, which generates two antipodal equilibria on $\S^2$. The same holds for $X_2$. The remaining trajectories of $X_1$ and $X_2$ are contained on invariant planes which intersect $\S^2$ on closed curves surrounding the equilibria, which corresponds to periodic orbits of those linear vector fields. In other words, $Z$ has two pairs of antipodes equilibria of center type on $\S^2$.

If $\s$ is a circle on $\S^2$ but not a great circle, then $Z$ has at least one real equilibrium point of center type and therefore it cannot be transitive. On the other hand, if $\s$ is a great circle, we have some situations to analyze. Indeed, if the equilibria of $X_1$ or $X_2$ lie outside $\Sigma$ then $Z$ has a real equilibria of center type so $Z$ cannot be transitive. Then assume the four equilibria of $Z$ lie on the great circle $\Sigma$ and notice that the points between the equilibria on $\Sigma$ may be of crossing or sliding/escaping type. If there is no sliding/escaping, then $Z$ has a pair of double elliptic tangency points. It is easy to see that trajectories around those tangency points behave like center equilibria, then again $Z$ cannot be transitive.

The last situation to study is the case of four distinct tangency points of $\Sigma$, generating four regions between the tangency points on $\Sigma$, being two crossing regions, one sliding region and an escaping one. In this case the trajectories associated to the sliding and escaping regions are invariant because their boundaries are formed by equilibria and therefore transitive cannot occur. This finishes the proof of the proposition.


	
	

\subsection{Proof of Theorem \ref{main:transitivity}}
\begin{proof}
We shall prove that $\Sigma^s$ and $\Sigma^e$ are non empty sets, hence we start with the following claim:

\begin{center}
\textit{Claim}: $\Sigma^e \neq \emptyset$.
\end{center}

The claim is proved by contradiction, then assume $\Sigma^e= \emptyset$. 

The set $\Sigma^s$ is either empty or nonempty, either case will lead to a contradiction. Indeed, notice that if the sliding set is empty, then we obtain a Filippov vector field without sliding and escaping region on the sphere. In particular it implies that the only way to a Filippov orbit to experience non-uniqueness of solution is at double tangency points. Hence consider a new vector field denoted by $\phi Z$, where $Z$ is the considered Filippov vector field and $\phi$ is a positive smooth function which is zero only at the tangency points which are finite by hypothesis. Now we have transformed the Filippov vector field $Z$ into a vector field with uniqueness of solutions.  

The previous procedure still leaves $\phi Z$ transitive. To see that, let $\gamma$ be a transitive orbit for $Z$. In order to be dense, $\gamma$ must visit these tangency points only a finite number of time and, at some point, it never returns to these tangency points. Hence, the forward orbit of $\gamma$ is dense and it is an orbit of $\phi Z$ (up to a reparametrization) since we only changed the dynamics on the tangency points. Therefore we obtain a continuous transitive vector field on the sphere, which is an absurd. Hence we have eliminated the empty slinding set. 

We still have to consider the case where the sliding set is nonempty. Due to transitivity, an orbit which enters the sliding region has to leave it eventually. Moreover, there exist only a finite number of exit points at the tangency ones. Since there exist no escaping region, such an orbit faces non-uniqueness only at tangency points which are finite. But then by transitivity that orbit must return the same sliding region and then it is a periodic orbit which is not dense. However it is an absurd which finishes the proof of the claim.

The claim implies that any Filippov vector field fulfilling the hypothesis of the theorem has nonempty escaping region. We now consider a new Filippov vector field replacing $X$ by $-X$. This is another way of inverting the orientation of orbits in time and for this new Filippov vector field we apply the claim we just proved, hence it has nonempty escaping region. But now the escaping region associated to $-X$ is the sliding region of the original Filippov vector field $X$. This proves item $a)$.

We now prove item $b)$. Consider a sliding and an escaping region, respectively, $I$ and $J$. Let $U$ and $V$ be two open sets close enough to $I$ and $J$, respectively, in such that any orbit starting on $U$ enters $I$ and any orbit passing through $V$ comes from $J$. By topological transitivity there exist an orbit that visits $U$ and $V$. Assume $\gamma$ is such an orbit with $\gamma(t_1) \in U$ and $\gamma(t_2) \in V$. We assume that $t_1 < t_2$, the other case can be treated by analogous arguments. Since $\gamma(t_1) \in U$ then before reaching $V$ it enters the sliding region and before the time $t_2$ it touches the escaping region, therefore the sliding and escaping regions $I$ and $J$ are connected. 

For the last part of item $b)$ observe that since sliding and escaping regions can be connected, we can consider an orbit that leave some escaping region and enter a sliding region. Notice that small perturbations of this orbit in the escaping region still generate an orbit entering the sliding region, hence we can create many different orbits connecting these regions.

\end{proof}

\subsection{Proof of Theorem \ref{main:no-transitive}}





\begin{proof}
The proof is done by contradiction. Consider $Z$ a PSVF with a finite number of tangency points and assume that $Z$ is robustly transitive. So there exist $U_Z$ a neighborhood of $Z$ such that every Filippov vector field in $U_Z$ is transitive.

We first notice that, from Theorem \ref{main:transitivity}, $Z$ has sliding and escaping regions, so take $p\in\Sigma^s$ and let $X$ and $Y$ be the adjacent vector fields separated by $\Sigma^s$. The positive sliding trajectory of $p$ eventually reaches a tangency point $T_A$ at the boundary of $\Sigma^s$, say that $X$ is tangent at $T_A$. We assume that $T_A$ is not a tangency point for $Y$ and that it is the common boundary between $\Sigma^s$ and a crossing region (see Figure \ref{fig3_thC}). Otherwise, we perturb $Z$ conveniently in order to obtain the described configuration. We also notice that, if $T_A$ is an invisible tangency point, then it is the $\omega-$limit set of every point nearby $\Sigma^s$, but this contradicts the transitivity of $Z$. Therefore $T_A$ is a visible tangency point of $Z$ and points on $\Sigma^s$ close to $T_A$ leave $\Sigma$ through this point. Moreover, this trajectory is unique (see Figure \ref{fig3_thC}).

\begin{figure}[h!]
	\includegraphics[scale=.7]{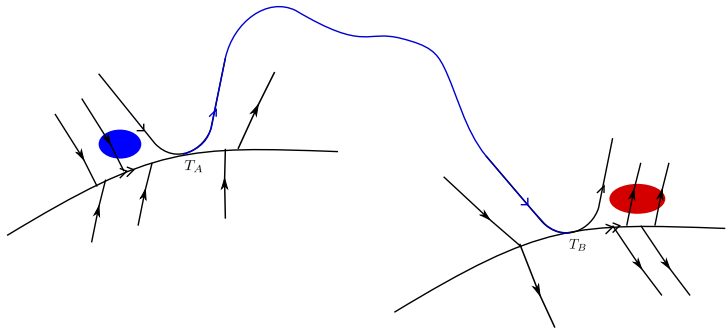}
	\caption{Connection between visible tangent points.}
	\label{fig3_thC}
\end{figure}

Analogously, there exist a visible tangency point $T_B$ located at the common boundary of some escape region $\Sigma^e$ and a crossing region which repels trajectories of the escaping vector field. In other words, $T_B$ is an entering point (the only one) to the escaping region $\Sigma^e$ (see Figure \ref{fig3_thC}).

The goal is to use the transitive property to connect $T_A$ to $T_B$ and then perturb $Z=Z_0$ in a suitable way to obtain a topologically transitive PSVF which does not connect tangency points (an absurd from Theorem \ref{main:transitivity}).


Effectively, let us call the escaping region associated to $T_B$ by $J$. Notice that if a point is in a neighborhood of $T_B$ and it is not in a trajectory which enters immediately in $J$ by the tangency $T_B$, then the time this point takes to enter $J$ is uniformly greater than some fixed time (this is because $J$ has a repelling behavior around it with the exception of the points which enters through $T_B$ the region $J$ itself). Let us call this number by $\alpha > 0$. Since $Z_0$ is transitive we know by Theorem \ref{main:transitivity} that we may connect through a Filippov orbit the tangency points $T_A$ and $T_B$.

Let $t_0 >0$ be the time for which an orbit from $T_A$ takes to reach $T_B$. Let us call this orbit by $\gamma_0$, that is, $\gamma_0$ is a Filippov orbit of $Z_0$ such that $\gamma_0(0)=T_A$, $\gamma_0(t_0)=T_B$ and $\gamma_0([0,t_0)) \cap \{T_B\}=\emptyset$.

Let $Z_1$ be a perturbation of $Z_0$ which has the following characteristics:
\begin{itemize}
\item[i)] $Z_1 \in U_Z$.
 \item[ii)] $Z_1$ coincides with $Z_0$ in $V_1^c$, where $V_1$ is an open ball which does not intersect $\Sigma$ (the switching manifold).
 \item[iii)] $\gamma_0(t)$ is an orbit of $Z_0$ and $Z_1$ as long as $\gamma_0[0,t] \subset V_1^c$.
 \item[iv)] let $\gamma_1$ be an orbit of  $Z_1$ which is a continuation of the orbit $\gamma_0$, then there exist a time $t_1 > t_0 + \alpha/2$ such that $\gamma_1([0,t_1))\cap \{T_B\} = \emptyset$, $\gamma_1(0) = T_A$ and $\gamma_1(t_1)=t_B$.
 \item[v)] $|Z_0-Z_1| < 1/2$.
\end{itemize}

The perturbation is done as follows. We call $J$ the escaping region associated to $T_B$ and let $\mathcal O$ be a small neighborhood of $J$ minus $J \cup T_B$ and minus the two connected segments that are inside this neighborhood which connects by a Filippov orbit with the tangency point $T_B$. In other words, $\mathcal O$ is a small region around $J$ for which the orbits have some sort of repelling behavior. And any point in $\mathcal O$ must take more than $\alpha$ to return to the region $J$, see Figure \ref{fig4_thC}.

\begin{figure}[h!]
	\includegraphics[scale=1]{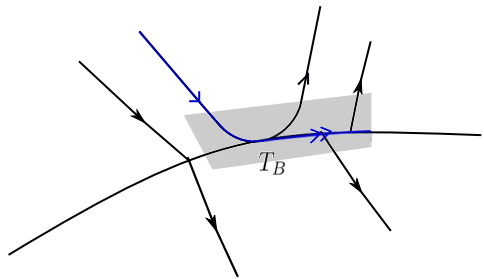}
	\caption{The neighborhood $\mathcal{O}$.}
	\label{fig4_thC}
\end{figure}


Let $\gamma_0(\xi_0)$ be some point very close to $T_B$ and $\xi_0 \in [0,t_0)$. By very close we mean that any ball around $\gamma_0(\xi_0)$ must intersect $\mathcal O$, in other words $\gamma_0(\xi_0)$ is so close to $T_B$ that by a perturbation small we can place it inside $\mathcal O$. Also, consider a ball $V_1$ around $\gamma_0(\xi_0)$ small enough such that its closure is inside the closure of $\mathcal O$, the closure of $V_1$ does not intersect $\Sigma$ and $\gamma_0$ enters $V_1$ in a time greater than $t_0-\alpha/2$.

We now perturb the vector field $Z_0$ inside $V_1$.  The perturbation will happen inside some compact set on $V_1$. Let us describe this perturbation. On $V_1$ we define $Z_1$ to be $Z_0 + W_1$, where $W_1$ is defined outside $V_1$ to be zero. Let us define it in $V_1$. We know that $Z_0(\gamma_0(\xi_0))$ is a nonzero vector. Consider $v_0$ a perpendicular vector to $Z_0(\gamma_0(\xi_0))$. Let $V_0$ be a smooth vector field around and define $W_1 := \phi V_0$ where $\phi$ is a smooth bump. We now consider the bump function to be sufficiently small in order that we can guarantee the conditions listed above.

We now proceed in a recursive way. Let $Z_n$ be a perturbation of $Z_{n-1}$ which has the following characteristics:
\begin{itemize}
\item[i)] $Z_n \in U_Z$.
 \item[ii)] $Z_n$ coincides with $Z_{n-1}$ in $V_n^c$, where $V_n$ is an open ball which does not intersect $\Sigma \cup V_1 \cup \ldots \cup V_{n-1}$.
 \item[iii)] $\gamma_{n-1}(t)$ is an orbit of $Z_{n-1}$ and $Z_n$ as long as $\gamma_{n-1}[0,t] \subset V_n^c$.
 \item[iv)] let $\gamma_n$ be an orbit of  $Z_n$ which is a continuation of the orbit $\gamma_{n-1}$, then there exist a time $t_n > t_{n-1} + \alpha/2$ such that $\gamma_n([0,t_n))\cap \{T_B\} = \emptyset$, $\gamma_n(0) = T_A$ and $\gamma_n(t_n)=T_B$.
 \item[v)] $|Z_{n}-Z_{n-1}| < 1/2^n$.
\end{itemize}

Note that $Z_n$ is uniformly converging to $ \widetilde Z$. Also $ \widetilde Z$ has never changed $Z$ on $\Sigma$. Since $ \widetilde Z$ is a transitive vector field it should connect $T_A$ and $T_B$, but the trajectory of $ \widetilde Z$ starting at $T_A$ which is the extension of $\gamma_0$ never touches $T_B$. This would give an absurd because we have to conect $T_A$ and $T_B$, but it turns out that the trajectory of $ \widetilde Z$ leaving $T_A$ and going to $T_B$ could be a different trajectory instead of the one which is an extension of $\gamma_0$ but the other possible way of leaving $T_A$. Hence on the proof above we would have to analyse at the same time both orbits. Thus, we get an absurd, that is, we have a transitive map which does not connect the tangency points $T_A$ and $T_B$.
\end{proof}

\subsection{Proof of Theorem \ref{main:general}}

\begin{proof}
    The proof is the same as the proof of Theorem \ref{main:no-transitive}. We observe that the only part of the proof for which $\mathbb S^2$ was really needed was to obtain that, for a transitive map on the sphere, one has necessarily sliding and escaping region. Hence that is the additional hypothesis on the theorem.
\end{proof}


\vspace{.3cm}

\noindent {\textbf{Acknowledgments.} This document is the result of the research projects funded by Pronex/ FAPEG/CNPq grant 2012 10 26 7000 803 and grant 2017 10 26 7000 508 (Euzébio), Capes grant 88881.068462/2014- 01 (Euzébio and Jucá), Universal/CNPq grant 420858/2016-4 (Euzébio), CAPES Programa de Demanda Social - DS (Jucá). R.V. was partially supported by National Council for Scientific and Technological Development – CNPq, Brazil and partially supported by FAPESP (Grants \#17/06463-3 and \# 16/22475-9).


\begin{thebibliography}{99}

\bibitem{BrouckePughSimic01} {\sc Broucke, Mireille E. and Pugh, Charles C. and Simi\'{c}, Slobodan N.}, {\it Structural stability of piecewise smooth systems}, Computational \& Applied Mathematics \textbf{20} (2001)

\bibitem{diBernardo-electrical-systems} {\sc M. di Bernardo, A. Colombo and E. Fossas}, {\it Two-fold singularity in nonsmooth electrical systems}, Proc. IEEE International Symposium on Circuits ans Systems, (2011), 2713--2716.

\bibitem{diBernardo-livro} {\sc M. di Bernardo, C.J. Budd, A.R. Champneys and P. Kowalczyk},
{\it Piecewise-smooth Dynamical Systems $-$ Theory and
	Applications}, Springer-Verlag (2008).

\bibitem{diBern-relay} {\sc M. di Bernardo, K.H.  Johansson, and F. Vasca}, {\it Self-Oscillations and Sliding in Relay Feedback Systems: Symmetry and Bifurcations}, Internat. J.  Bif. Chaos, vol. 11, (2001), pp. 1121-1140.

\bibitem{Brogliato} {\sc B. Brogliato}, {\it Nonsmooth Mechanics: Models, Dynamics and Control}, Springer- Verlag, New York, 1999.

\bibitem{BCE} {\sc C.A. Buzzi, T. Carvalho and R.D. Euz\'{e}bio},
{\textit On Poincaré-Bendixson Theorem and non-trivial minimal sets in planar nonsmooth vector fields}, Publicacions Mathematiques, Vol.~62, (2018), 113--131.
	
	
\bibitem{BCE-ETDS} {\sc C.A. Buzzi, T. Carvalho and R.D. Euz\'{e}bio},
{\it Chaotic Planar Piecewise Smooth Vector Fields With Non Trivial Minimal Sets}, Ergodic Theory of Dynamical Systems, Vol.~36, (2016), 458--469. 

\bibitem{Carvalho-LFernando} {\sc Carvalho, Tiago and Gon\c{c}alves, Luiz Fernando}, {\it Combing the {H}airy {B}all {U}sing a {V}ector {F}ield {W}ithout {E}quilibria}, {J. Dyn. Control Syst.}, v. 26, pp. 233-242 (2020).

\bibitem{Rossa} {\sc F. Dercole and F.D. Rossa}, {\it Generic and Generalized Boundary Operating Points in Piecewise-Linear (discontinuous) Control Systems}, In 51st IEEE Conference on Decision and Control, 10-13 Dec. 2012, Maui, HI, USA. Pages:7714-7719.

\bibitem{devaney} {\sc Devaney, Robert L}. {\it An introduction to chaotic dynamical systems}, Studies in Nonlinearity. Westview Press, Boulder, CO, 2003.

\bibitem{Dixon} {\sc D.D. Dixon}, {\it Piecewise Deterministic Dynamics from the Application of Noise to Singular Equation of Motion}, J. Phys A: Math. Gen. 28, (1995), 5539--5551.

\bibitem{EV} {\sc R.D. Euzébio and R. Varão},
{\it Topological transitivity imply chaos for two--dimensional Filippov systems}, submitted, 2019.

\bibitem{Fi} {\sc A.F. Filippov},
{\it Differential Equations with Discontinuous Righthand Sides},
Mathematics and its Applications (Soviet Series), Kluwer Academic
Publishers-Dordrecht, 1988.

\bibitem{Genema} {\sc S. Genena, D.J. Pagano and P.  Kowalczik}, {\it Hosm Control of Stick-Slip Oscillations in Oil Well Drill-Strings}, In Proceedings of the European Control Conference 2007 - ECC07, Kos, Greece, July.

\bibitem{Marcel} {\sc M. Guardia, T.M. Seara and M.A.
Teixeira}, {\it Generic bifurcations of low codimension of planar
Filippov Systems}, Journal of Differential Equations \textbf{250}
(2011) 1967--2023.

\bibitem{Jac-To} {\sc A. Jacquemard and D.J. Tonon}, {\sc  Coupled systems of non-smooth differential equations}, Bulletin des Sciences Math\'{e}matiques, vol. 136, (2012), pp. 239-255.

\bibitem{Lopez-Lopez04} {\sc Jim\'{e}nez L\'{o}pez, V\'{\i}ctor and Soler L\'{o}pez, Gabriel}, {\sc Transitive flows on manifolds},
Rev. Mat. Iberoamericana, vol. 20, (2004), pp. 107-130.

\bibitem{Kousaka} {\sc T. Kousaka, T.  Kido, T.  Ueta, H. Kawakami and M. Abe}, {\it Analysis of Border- Collision Bifurcation in a Simple Circuit}, Proceedings of the International Symposium on Circuits and Systems, II-481-II-484, (2000).

\bibitem{Leine} {\sc R. Leine and H.  Nijmeijer}, {\it Dynamics and Bifurcations of Non-Smooth Mechanical Systems}, Lecture Notes in Applied and Computational Mechanics, vol. 18, Berlin Heidelberg New-York, Springer-Verlag, (2004).

\bibitem{Smith-Russel88} {\sc Smith, Russell A. and Thomas, E. S.},
{\it Transitive flows on two-dimensional manifolds}, J. London Math. Soc. \textbf{(2)}, 569-576, (1988).

\bibitem{Utkin1992} {\sc Utkin, Vadim I.}, {\it Sliding modes in control and optimization}, Communications and Control Engineering Series, Springer-Verlag,  1992

\bibitem{Barbashin1970} {\sc Barbashin, E. A.}, {\it Introduction to the theory of stability},  {Wolters-Noordhoff Publishing, Groningen}, 1970.
\end{thebibliography}
\end{document}